\documentclass[a4paper,reqno]{amsart}
\usepackage[english]{babel}
\usepackage{amsmath,amsthm,amssymb,amsfonts,xcolor}

\usepackage{enumitem}
\setlist[enumerate]{leftmargin=*}
\setlist[itemize]{leftmargin=*}

\usepackage{mathrsfs}
\usepackage{hyperref}

\numberwithin{equation}{section}

\newtheorem{theorem}{Theorem}[section]
\newtheorem{lemma}[theorem]{Lemma}

\theoremstyle{definition}

\newcommand{\TT}{T}
\newcommand{\ee}{\mathrm{e}}
\newcommand{\dd}{\mathrm{d}}
\newcommand{\CC}{\mathbb{C}}
\newcommand{\ZZ}{\mathbb{Z}}

\newcommand{\RR}{\mathcal{R}}
\newcommand{\R}{\mathbb{R}}
\newcommand{\LL}{\mathcal{L}}
\newcommand{\var}{\varepsilon}

\newcommand{\hardy}{H_{\mathrm{at}}^1(\mu)}
\newcommand{\defeq}{\mathrel{:=}}     
\newcommand{\eqdef}{\mathrel{=:}}

\title[Endpoint behaviour of Riesz transforms on homogeneous trees]{Endpoint behaviour of vertical and horizontal Riesz transforms for the flow Laplacian on homogeneous trees}

\author{Jie Liu}
\address[J. Liu]{School of Mathematics and Statistics \\ Northwestern Polytechnical University \\ Xi'an 710129 \\ China \\ and Dipartimento di Scienze Matematiche \\ Politecnico di Torino \\ Corso Duca degli Abruzzi 24 \\ 10129 Torino \\ Italy}
\email{jay2000@mail.nwpu.edu.cn}
\email{jie.liu@polito.it}


	
\thanks{The author gratefully acknowledges the financial support of the China Scholarship Council (Grant No.\ 202406290140).}

\begin{document}

\begin{abstract}
We study endpoint estimates of Riesz transforms on a homogeneous tree equipped with its canonical flow measure $\mu$. We prove that the known boundedness from $\hardy$ to $L^1(\mu)$ of the vertical Riesz transform can not be strengthened to boundedness on $\hardy$ itself. Using an explicit family of atoms and a logarithmic $\mathrm{BMO}$ test function, we obtain logarithmically growing lower bounds for the Hardy norms of their images. We also prove that the adjoint vertical Riesz transform is not of weak type $(1,1)$. For the scalar horizontal transforms defined by nonzero bounded coefficients with zero sum on each child set, we establish $\hardy\to L^1(\mu)$ boundedness and show that failure to preserve cancellation on atoms prevents boundedness on $\hardy$. We then introduce a horizontal derivative based on deviations from sibling averages. Its associated Riesz transform retains the $L^p(\mu)$ bounds for $1<p<\infty$ and the weak type $(1,1)$ estimate, and is bounded on $\hardy$ itself. These mapping properties parallel those of horizontal Riesz transforms in the continuous setting, supporting our construction as a natural discrete analogue. These results distinguish vertical and horizontal behaviour at the endpoint and identify cancellation on sibling sets as a mechanism for Hardy-space boundedness.
\end{abstract}

\maketitle

\section{Introduction}
Riesz transforms play a central role in the real-variable theory of Hardy spaces. In the Euclidean setting, the first-order Riesz transforms $R_j=\partial_j(-\Delta)^{-1/2}$ are of weak type $(1,1)$, bounded on $H^1(\R^n)$, and characterize $H^1(\R^n)$ through the norm equivalence
\[
 \|f\|_{H^1(\R^n)}\simeq \|f\|_{L^1(\R^n)}+\sum_{j=1}^n\|R_jf\|_{L^1(\R^n)},
\]
and are themselves bounded on $H^1(\R^n)$; see \cite{FS72} and \cite[Chapter~III]{S93}. Moreover, $R_j^*=-R_j$, so the transforms and their adjoints have the same mapping properties. On spaces of exponential growth, these conclusions can separate: boundedness from an atomic Hardy space to $L^1$ need not imply boundedness on the Hardy space itself, and a weak type estimate for a Riesz transform need not hold for its adjoint. This paper studies these distinctions on homogeneous trees with their canonical flow measure, and examines how the choice of horizontal differentiation affects the Hardy-space endpoint.

Let $\TT=\TT_{q+1}$, $q\ge2$, be a homogeneous tree oriented towards a fixed end. For each vertex $x$, we denote by $p(x)$ its parent and by $s(x)$ its set of $q$ children. Fix a reference vertex $o$ and normalize the level function $\ell$ so that $\ell(o)=0$ and $\ell(p(x))=\ell(x)+1$. The canonical flow measure
\[
 \mu(x)=q^{\ell(x)}
\]
then satisfies $\mu(o)=1$ and $\mu(x)=\sum_{y\in s(x)}\mu(y)$. With the graph distance, this measured tree has exponential volume growth and is nondoubling. 

The classical atomic Hardy-space theory of Coifman and Weiss \cite{CW77} is formulated on spaces of homogeneous type and therefore does not apply directly to the measured tree considered here. To deal with this geometry, Hebisch and Steger \cite{HS03} developed a Calder\'on--Zygmund theory based on suitable trapezoids rather than metric balls. Building on this framework, Arditti, Tabacco and Vallarino \cite{ATV20,ATV21} introduced the corresponding atomic Hardy and $\mathrm{BMO}$ spaces and established their duality and interpolation properties.

More recently, Conde Alonso, Martini and Vallarino \cite{CMV26} developed a general nondoubling Calder\'on--Zygmund theory based on families of admissible sets equipped with recursive partitions. Their framework provides atomic Hardy and $\mathrm{BMO}$ spaces, duality and interpolation results, as well as endpoint estimates for singular integral operators satisfying suitable kernel conditions. In particular, it applies to arbitrary flow trees without requiring their flow measures to be locally or globally doubling; see \cite[Section~4]{CMV26}. Throughout this paper, $\hardy$ denotes the atomic Hardy space associated with admissible trapezoids on the tree under consideration.

Let \(\nabla\) denote the flow gradient, defined by
\begin{equation}\label{Def-nabla}
\nabla f(x)=f(x)-f(p(x)),
\end{equation}
and let \(\LL=\nabla^*\nabla/2\) be the flow Laplacian.
The associated vertical Riesz transform is defined by
\begin{equation}\label{Def-verRiesz}
\RR=\nabla\LL^{-1/2}.
\end{equation}
This Riesz transform has been extensively studied. Its weak type $(1,1)$ estimate and $L^p(\mu)$ boundedness for $1<p\le 2$  follow from \cite[Theorem~2.3]{HS03}. The full range of $L^p(\mu)$ boundedness, $1<p<\infty$, was subsequently established in \cite[Theorem~1.1]{LMSTV23}. At the Hardy-space endpoint, the boundedness $\RR:\hardy\to L^1(\mu)$ is also known; see \cite[Theorem~2.3]{HS03}, \cite[Theorem~3]{ATV20}, or \cite[Theorem~1.2~(1)]{S22}.

These results have since been extended to more general flow trees. Levi, Santagati, Tabacco and Vallarino \cite{LSTV23} developed a Calder\'on--Zygmund theory and a Hardy--$\mathrm{BMO}$ framework for locally doubling flow measures. More recently, Martini, Santagati, Tabacco and Vallarino \cite[Theorem~1.2]{MSTV25} established $L^p$ boundedness for the full range $1<p<\infty$ on locally finite flow trees, together with a sharp Mihlin--H\"ormander spectral multiplier theorem. When the flow measure is additionally assumed to be locally doubling, they also obtained weak type $(1,1)$ and $\hardy\to L^1(\mu)$ estimates. Analogous endpoint estimates are established in \cite{CMV26} even without the local doubling assumption, with an appropriate definition of the atomic Hardy space.

Continuous models for related questions arise on solvable Lie groups of exponential growth. On the $ax+b$ group $G_n=\R^n\rtimes\R^+$, the distinguished Laplacian and its first-order Riesz transforms are
\[
 \LL_G=-\sum_{j=0}^n X_j^2,\,\,\,\RR_{G,j}=X_j\mathcal L_G^{-1/2},
\]
where $X_0=a\partial_a$ is the vertical field and $X_j=a\partial_{x_j}$, $1\le j\le n$, are the horizontal fields. The dilation direction and horizontal level structure have discrete counterparts in the orientation and horocycles of a homogeneous
flow tree; see \cite{LMSTV23}. Martini \cite[Theorem~1.1]{M23} proved $L^p$ boundedness of $\RR_{G,j}$ for every $1<p<\infty$ and weak type $(1,1)$ estimate from \cite{HS03}. More generally, Liu and Martini \cite{LM26} obtained the same range of $L^p$ boundedness, as well as weak type $(1,1)$ estimates, for the first-order Riesz transforms associated with the distinguished Laplacian on Damek--Ricci spaces, with respect to right Haar measure. 

The Hardy-space framework on these groups was developed by Vallarino \cite{V09}, who established Hardy--BMO duality and interpolation results adapted to their nondoubling geometry. On $G_2$, Sj\"ogren and Vallarino \cite{SV08} proved that all three first-order Riesz transforms are bounded from the atomic Hardy space to $L^1$, whereas none of their adjoints enjoys this property. They subsequently showed that these Riesz transforms do not characterize the associated atomic Hardy space; see \cite[Theorem~1.3]{SV25}. They also proved that both horizontal transforms are bounded on this Hardy space \cite[Theorem~1.4]{SV25}, leaving the corresponding question for the vertical transform open.

On homogeneous flow trees, Santagati \cite[Theorem~1.2]{S22} likewise proved that the Riesz-transform characterization fails:
\[
 \hardy\subsetneq H^1_{\RR}(\mu)\defeq\{f\in L^1(\mu):\RR f\in L^1(\mu)\}.
\]
This shows that integrability of $f$ and $\RR f$ does not characterize membership in the atomic Hardy space. It leaves a different question: whether $\RR$ preserves that space when the input already belongs to $\hardy$. Our first main result answers this question negatively.
\begin{theorem}\label{Thm-mainver}
There are $(1,\infty)$-atoms $a_h$, indexed by integers $h\ge2$, and constants $c_q,C_q>0$ such that
\[
 \|\RR a_h\|_{\hardy}\ge c_q\log h-C_q.
\]
Here the norm is interpreted as $+\infty$ when $\RR a_h\notin\hardy$. Consequently, $\RR$ does not extend to a bounded operator from $\hardy$ to itself.
\end{theorem}

Thus the failure of boundedness persists despite both the known $\hardy\to L^1(\mu)$ estimate and preservation of the zero integral condition: these properties do not control the atomic Hardy norm.

The proof relates the growth of these norms to the geometry of the transform away from the support of the atoms. We choose $a_h$ on two consecutive blocks of descendant layers below $o$ and test $\RR a_h$ against bounded truncations of
\[
 b(x)=\log\bigl(1+d(x,\gamma)\bigr),\,\,\,\gamma=\{p^k(o):k\ge0\}.
\]
These truncations have uniformly bounded $\mathrm{BMO}$ seminorms. Outside the descendant subtree of $o$, the transform is positive on the ancestor ray $\gamma$ and negative on its side branches. The test function vanishes on the ray, while its logarithmic growth detects the branch contribution at distances comparable to $h$. A summation-by-parts argument bounds the contribution from the descendant subtree uniformly in $h$. The resulting pairings yield a logarithmic lower bound for the Hardy norm through Hardy--$\mathrm{BMO}$ duality.

We next turn to the endpoint behaviour of the adjoint $\RR^*$. It was already known that $\RR^*$ is not bounded from $\hardy$ to $L^1(\mu)$, as follows from
\cite[Proposition~4.4]{LMSTV23} and Hardy--$\mathrm{BMO}$ duality in \cite[Theorem~13]{ATV21} or \cite[Theorem~4.10]{LSTV23}. Analogous failures of $H^1\to L^1$ boundedness
are known on $ax+b$ groups \cite[Proposition~7.1]{M23} and Damek--Ricci spaces \cite{BL26}.

At the weak type $(1,1)$ endpoint, a distinction arises between horizontal and vertical Riesz transforms. Levi, Martini, Santagati, Tabacco and Vallarino \cite[Theorem~5.5]{LMSTV23} proved that the adjoints of the horizontal Riesz transforms on homogeneous trees are of weak type $(1,1)$, while leaving the corresponding question for $\RR^*$ open \cite[Remark~4.2]{LMSTV23}.

Our second main result answers this question negatively: $\RR^*$ is not of weak type $(1,1)$. This establishes a distinction between the horizontal and vertical adjoints and provides a discrete counterpart to Hormozi's result \cite[Theorem~1.1]{H26}.

\begin{theorem}\label{Thm-mainadjver}
The adjoint vertical Riesz transform $\RR^*$ is not of weak type $(1,1)$.
\end{theorem}
Motivated by the continuous counterexample in the proof of \cite[Theorem~1.1]{H26}, we construct nonnegative, finitely supported functions of unit $L^1$ norm whose mass is distributed sparsely over many descendant layers. The explicit formula for $\RR-\RR^*$ separates descendant averages from ancestor values; see \cite[Eq.~(3.4) and (4.4)]{LMSTV23}. Removing the latter contribution on a sufficiently large set produces logarithmically growing weak $L^1$ norms. The known weak type $(1,1)$ bound for $\RR$ then yields the conclusion for $\RR^*$.

We next turn to horizontal differentiation. On $ax+b$ groups, the semidirect-product structure distinguishes the dilation direction from the translation directions,
giving rise to vertical and horizontal left-invariant derivatives. On a tree, the chosen orientation defines the parent derivative in \eqref{Def-nabla}, but it does not select horizontal coordinate directions: vertices on the same horocycle are not joined by edges. A horizontal derivative must therefore specify how values within a level are compared.

Motivated by the continuous theory, Levi, Martini, Santagati, Tabacco and Vallarino introduced in \cite[Section~5]{LMSTV23} a family of scalar horizontal Riesz transforms $\RR_\var$. Specifically speaking, for a nonzero $\var\in L^\infty(\mu)$ satisfying
\begin{equation}\label{Def-var}
 \sum_{y\in s(x)}\var(y)=0,\,\,\, x\in\TT,
\end{equation}
the scalar horizontal Riesz transform is defined by
\begin{equation}\label{Def-scalhorRiesz}
 \nabla_\var f(x)=\frac1q\sum_{y\in s(x)}\var(y)f(y),\,\,\,\RR_\var=\nabla_\var\LL^{-1/2}.
\end{equation}
The children of a vertex lie on the same horocycle, so this derivative compares values within a horizontal level. These transforms are bounded on $L^p(\mu)$ for $1<p<\infty$
and are of weak type $(1,1)$; their adjoints are also of weak type $(1,1)$; see \cite[Section~5]{LMSTV23}. Nevertheless, their behaviour on the atomic Hardy space is different from that of the continuous horizontal transforms on $G_2$.  \begin{theorem}\label{Thm-varH1}
For the above nonzero \(\var\), \(\RR_\var\) is bounded from \(\hardy\) to \(L^1(\mu)\), but not bounded from \(\hardy\) to itself.
\end{theorem}
Indeed, for every nonzero $\var$ as above, there is a $(1,\infty)$-atom $a$ such that
\[
 \sum_{x\in\TT}\RR_\var a(x)\mu(x)\ne0.
\]
Thus, although the coefficient condition makes $\nabla_\var$ annihilate constants, it does not ensure that $\RR_\var$ preserves the zero integral condition of the atomic space.

This obstruction motivates a horizontal derivative whose output has cancellation on each sibling set \(s(p(x))\). We define
\[
\nabla_h f(x)=\frac1q\left(f(x)-\frac1{q-1}\sum_{z\in b(x)}f(z)\right),\,\,\,\RR_h =\nabla_h\LL^{-1/2},
\]
where \(b(x)=s(p(x))\setminus\{x\}\) is the set containing all the siblings of \(x\). 
For later use, we rewrite this operator in a form which makes the cancellation explicit. Set \(u=\LL^{-1/2}f\), then 
\[
\RR_h f(x)=\nabla_h u(x)=\frac1q\left(u(x)-\frac1{q-1}\sum_{z\in b(x)}u(z)\right),
\]
and thus,
\[
\sum_{x\in s(y)}\RR_h f(x)=0.
\]
Thus the horizontal Riesz transform has a built-in cancellation across the \(q\) children of each vertex. The associated Riesz transform retains the $L^p$, weak type
$(1,1)$, and $\hardy\to L^1(\mu)$ bounds of the scalar construction, and satisfies the stronger Hardy-space estimate below.
\begin{theorem}\label{Thm-mainhor} 
The horizontal Riesz transform \(\RR_h\) is bounded from \(\hardy\) to itself. 
\end{theorem}
The comparison of the two horizontal constructions highlights the role of cancellation on sibling sets. The symmetric construction $\RR_h$ enjoys this local cancellation and satisfies the positive Hardy-space endpoint estimate, providing a discrete analogue of the corresponding result on $\R^2\rtimes\R^+$ in \cite[\text{Theorem}~1.4]{SV25}, whereas the scalar construction does not. This comparison concerns their mapping properties and does not require choosing a preferred horizontal coordinate on the tree.

Together with the previously known estimates, our results complete the classification of the three endpoint mapping properties listed in Table~\ref{Tab-endpoints} for the operators considered here. All four operators are bounded on $L^p(\mu)$ for $1<p<\infty$.

\begin{table}[htbp]
\centering
\small
\setlength{\tabcolsep}{3pt}
\begin{tabular}{lccc}
\hline
Operator
& $L^1(\mu)\to L^{1,\infty}(\mu)$
& $\hardy\to L^1(\mu)$
& $\hardy\to\hardy$ \\
\hline
$\RR$       & Yes, \cite[Theorem~2.3]{HS03} & Yes, \cite[Theorem~3]{ATV20} & No, Theorem~\ref{Thm-mainver} \\
$\RR^*$     & No, Theorem~\ref{Thm-mainadjver} & No, \cite[Proposition~4.4]{LMSTV23} & No, \cite[Proposition~4.4]{LMSTV23} \\
& &and \cite[Theorem~13]{ATV21} &and \cite[Theorem~13]{ATV21} \\
$\RR_\var$  & Yes, \cite[Section~5]{LMSTV23} & Yes, Theorem~\ref{Thm-varH1} & No, Theorem~\ref{Thm-varH1} \\
$\RR_h$     & Yes, Theorem~\ref{Thm-Lpweakttype11q} & Yes, Theorem~\ref{Thm-Lpweakttype11q} & Yes, Theorem~\ref{Thm-mainhor} \\
\hline
\end{tabular}
\caption{Endpoint boundedness for Riesz transforms.}\label{Tab-endpoints}
\end{table}

The paper is organized as follows. Section~\ref{Sec-pri} collects the geometric and analytic preliminaries. The vertical Hardy-space counterexample is proved in Section~\ref{Sec-ver}. The failure of weak type $(1,1)$ for the adjoint vertical transform is proved in Section~\ref{Sec-adjver}. Section~\ref{Sec-hor} treats the
scalar horizontal transform and the transform defined through sibling averages.

\section{Preliminaries}\label{Sec-pri}

\subsection{Geometry and the canonical flow measure}

Let $\TT=\TT_{q+1}$, $q\ge2$, be a homogeneous tree of degree $q+1$, equipped with the graph distance $d$. We orient $\TT$ towards a fixed end. For each vertex $x$, we denote by $p(x)$ its parent, namely the unique neighbour in the direction of the fixed end, and by $s(x)$ its set of $q$ children. We write $y\le x$ if $x=p^k(y)$ for some integer $k\ge0$.

Fix a reference vertex $o\in\TT$. The level function $\ell:\TT\to\ZZ$ is normalized by
\[
 \ell(o)=0,\,\,\,\ell(p(x))=\ell(x)+1.
\]
We equip $\TT$ with the canonical flow measure
\[
 \mu(x)=q^{\ell(x)},\,\,\, x\in\TT,
\]
and write $\mu(E)=\sum_{x\in E}\mu(x)$ for $E\subseteq\TT$. In particular, $\mu(o)=1$, and
\[
 \mu(x)=\sum_{y\in s(x)}\mu(y).
\]

For a vertex $x$ and an integer $s\ge0$, define the descendant
layer
\begin{equation}\label{Def-descendantlayer}
 S_s(x)=\{y\in\TT:y\le x,\ d(x,y)=s\}.
\end{equation}
This layer consists of $q^s$ vertices, each of measure
$q^{-s}\mu(x)$. Consequently,
\[
 \mu(S_s(x))=\mu(x),\,\,\,x\in\TT,\ s\ge0.
\]

\subsection{Atomic Hardy and BMO spaces} 
A subset $R\subseteq\TT$ is called an admissible trapezoid if it is either a singleton or a set of the form
\[
 R=R_{h'}^{h''}(y_R)=\{y\in\TT:y\le y_R,\ h'\le d(y,y_R)<h''\},
\]
where $y_R\in\TT$ and $h',h''$ are positive integers satisfying
\[
 2\le\frac{h''}{h'}\le12.
\]
We denote by $\mathcal F$ the family of all admissible trapezoids. Since each descendant layer below $y_R$ has measure $\mu(y_R)$, every nonsingleton trapezoid satisfies
\[
 \mu(R)=(h''-h')\mu(y_R).
\]

A function $a:\TT\to\CC$ is a $(1,\infty)$-atom if there exists $R\in\mathcal F$ such that
\[
 \operatorname{supp}a\subseteq R,\,\,\,\|a\|_{L^\infty(\mu)}\le\frac1{\mu(R)},\,\,\,
 \sum_{x\in R}a(x)\mu(x)=0.
\]
The atomic Hardy space $\hardy$ consists of all $f\in L^1(\mu)$ admitting a representation
\[
 f=\sum_j\lambda_j a_j\,\,\text{in }L^1(\mu),\,\,\,\sum_j|\lambda_j|<\infty,
\]
where each $a_j$ is a $(1,\infty)$-atom. Its norm is
\[
 \|f\|_{\hardy}=\inf\sum_j|\lambda_j|,
\]
where the infimum is taken over all atomic decompositions of $f$. In particular,
\[
 \|f\|_{L^1(\mu)}\le\|f\|_{\hardy},\,\,\,\sum_{x\in\TT}f(x)\mu(x)=0,\,\,\, f\in\hardy.
\]
We refer to \cite[Section~2.2]{S22} or \cite[Section~4.2]{LSTV23} for further details.

For a function $b:\TT\to\CC$ and $R\in\mathcal F$, write
\[
 b_R=\frac1{\mu(R)}\sum_{x\in R}b(x)\mu(x).
\]
The space $\mathrm{BMO}(\mu)$ consists of functions with finite seminorm
\[
 \|b\|_{\mathrm{BMO}(\mu)}=\sup_{R\in\mathcal F}\frac1{\mu(R)}\sum_{x\in R}|b(x)-b_R|\mu(x).
\]
Modulo constants, this space is identified with the dual of $\hardy$; see \cite[Theorem~13]{ATV21} or \cite[Theorem~4.10]{LSTV23}. 
For $f\in\hardy$ and $b\in\mathrm{BMO}(\mu)\cap L^\infty(\mu)$, we shall use the pairing estimate
\[
 \left|\sum_{x\in\TT}f(x)b(x)\mu(x)\right|\le\|f\|_{\hardy}\|b\|_{\mathrm{BMO}(\mu)}.
\]
Indeed, for an atom $a$ supported on $R$, cancellation gives
\[
 \left|\sum_{x\in\TT}a(x)b(x)\mu(x)\right|=\left|\sum_{x\in R}a(x)(b(x)-b_R)\mu(x)\right|\le\|b\|_{\mathrm{BMO}(\mu)},
\]
and the general estimate follows by atomic decomposition. The boundedness assumption on $b$ ensures that the pairing is an absolutely convergent sum.

We record two elementary properties of $\mathrm{BMO}(\mu)$ that will be used below. First, for every $R\in\mathcal F$ and every constant $c_R\in\CC$,
\begin{equation}\label{Ieq-bmocomparison}
\frac1{\mu(R)}\sum_{x\in R}|b(x)-b_R|\mu(x)\le\frac2{\mu(R)}\sum_{x\in R}|b(x)-c_R|\mu(x).
\end{equation}
Thus, to prove that $b\in\mathrm{BMO}(\mu)$, it suffices to choose constants $c_R$ such that the average deviations from $c_R$ are bounded uniformly over $R\in\mathcal F$; see, for instance, \cite[Proposition~3.2~(ii)]{CMV26}.

Second, if $b\ge0$ belongs to $\mathrm{BMO}(\mu)$, then the bounded truncations
\[
b_M(x)=\min\{b(x),M\},\,\,\, M>0,
\]
satisfy
\begin{equation}\label{Ieq:bmotruncation}
\sup_{M>0}\|b_M\|_{\mathrm{BMO}(\mu)}\le 2\|b\|_{\mathrm{BMO}(\mu)};
\end{equation}
see, for instance, \cite[Proposition~3.6~(i)]{CMV26}.

\subsection{The flow Laplacian and its potential kernel}

The flow Laplacian on $(\TT,\mu)$ is given by
\[
 \LL f(x)=f(x)-\frac12f(p(x))-\frac1{2q}\sum_{y\in s(x)}f(y);
\]
see, for instance, \cite[Section~2.3]{S22}. For the vertical derivative
\begin{equation*}
 \nabla f(x)=f(x)-f(p(x)),
\end{equation*}
its adjoint with respect to $L^2(\mu)$ is
\[
 \nabla^*f(x)=f(x)-\frac1q\sum_{y\in s(x)}f(y).
\]
Thus $\LL=\frac12\nabla^*\nabla$, and $\LL$ is nonnegative
and self-adjoint on $L^2(\mu)$. The associated vertical
Riesz transform is
\[
 \RR=\nabla\LL^{-1/2}.
\]

We recall the kernel formula for $\LL^{-1/2}$ from
\cite[Proposition~2.3 and Corollary~2.5]{LMSTV23}.
For finitely supported $f$,
\[
 \LL^{-1/2}f(x)
 =\sum_{y\in\TT}K_{\LL^{-1/2}}(x,y)f(y)\mu(y),
\]
where
\begin{equation}\label{Kernel-LL1}
 K_{\LL^{-1/2}}(x,y)
 =q^{-\ell(x)/2}G(d(x,y))q^{-\ell(y)/2}.
\end{equation}
For every nonnegative integer $n$,
\begin{equation}\label{For-G}
 G(n)=\sum_{k=0}^{\infty}q^{-(n+2k)/2}\widetilde k_{\ZZ}(n+2k+1),\,\,\,
 \widetilde k_{\ZZ}(n)=\frac{2\sqrt2}{\pi}\frac{n}{n^2-1/4},
\end{equation}
where $\widetilde k_{\ZZ}$ is the convolution kernel of the skew-symmetric part of the Riesz transform on $\ZZ$. 
Subtracting the series for $G(n+2)$ from that for $G(n)$ gives
\begin{equation}\label{For-Gnn+2}
 q^{n/2}\bigl(G(n)-G(n+2)\bigr)=\widetilde k_{\ZZ}(n+1),\,\,\, n\ge0.
\end{equation}

\subsection{Shift operators and the scalar horizontal derivative}
Define the shift operator $\Sigma$ and its adjoint on $L^2(\mu)$ by
\begin{equation}\label{Def-shift}
\Sigma f(x)=f(p(x)),\,\,\,\Sigma^*f(x)=\frac1q\sum_{y\in s(x)}f(y).
\end{equation}
For every $1\le p\le\infty$, we have
\begin{equation}\label{shift}
\|\Sigma f\|_{L^p(\mu)}=\|f\|_{L^p(\mu)},\,\,\,\|\Sigma^*f\|_{L^p(\mu)}\le\|f\|_{L^p(\mu)}.
\end{equation}
Moreover,
\begin{equation}\label{shiftad}
\|\Sigma f\|_{L^{1,\infty}(\mu)}=\|f\|_{L^{1,\infty}(\mu)},\,\,\,\|\Sigma^*f\|_{L^{1,\infty}(\mu)}\le q\|f\|_{L^{1,\infty}(\mu)}.
\end{equation}
See \cite[Proposition~2.1]{LMSTV23} and, for the more general setting of flow trees, \cite[Section~2]{MSTV25}.

Let $\var\in L^\infty(\mu)$ satisfy
\begin{equation}\label{eq:var-cancellation}
\Sigma^*\var=0,\,\,\,\text{equivalently}\,\,\,\sum_{y\in s(x)}\var(y)=0
\,\,\,\text{for every }x\in\TT.
\end{equation}
The scalar horizontal derivative associated with $\var$ is defined by
\begin{equation}\label{Def-horizontalscalar}
\nabla_\var f(x)=\Sigma^*(\var f)(x)=\frac1q\sum_{y\in s(x)}\var(y)f(y).
\end{equation}
The associated scalar horizontal Riesz transform is defined by \(\RR_\var=\nabla_\var \LL^{-1/2}\).
The estimates \eqref{shift} and \eqref{shiftad} yield \begin{equation}\label{Est-nablavar}
\|\nabla_\var f\|_{L^p(\mu)}\le\|\var\|_\infty\|f\|_{L^p(\mu)},\,\,\, 1\le p\le\infty,
\end{equation}
and
\[
\|\nabla_\var f\|_{L^{1,\infty}(\mu)}\le q\|\var\|_\infty\|f\|_{L^{1,\infty}(\mu)};
\]
see also \cite[Proposition~5.1]{LMSTV23}. The cancellation condition \eqref{eq:var-cancellation} implies that $\nabla_\var\Sigma=0$. Recalling from \eqref{Def-nabla} that 
\begin{equation}\label{Rel-nablaSigma}
    \nabla=I-\Sigma,
\end{equation}
we obtain 
\begin{equation}\label{Rel-derscalarvert}
\nabla_\var=\nabla_\var\nabla.
\end{equation}
Thus the scalar horizontal Riesz transform, defined by $\RR_\var=\nabla_\var\LL^{-1/2}$, admits the factorisation
\begin{equation}\label{Rel-scalarvert}
\RR_\var=\nabla_\var\nabla\LL^{-1/2}=\nabla_\var\RR.
\end{equation}
Together with the preceding bounds for $\nabla_\var$, this identity shows that $\RR_\var$ inherits the known $L^p(\mu)$ boundedness for $1<p<\infty$, weak type $(1,1)$ estimate, and boundedness from $\hardy$ to $L^1(\mu)$
of the vertical Riesz transform $\RR$; see \cite[Section~5]{LMSTV23}.

\section{Unboundedness of the vertical Riesz transform on \texorpdfstring{$\hardy$}{H1}}\label{Sec-ver}

In this section, we will prove Theorem~\ref{Thm-mainver}: $\hardy\to\hardy$ unboundedness of $\RR$. The proof combines an explicit family of atoms with a logarithmic $\mathrm{BMO}$ test function adapted to the geometry of the tree.
\begin{proof}[Proof of Theorem~\ref{Thm-mainver}]
\textbf{Step 1: construction of atoms $a_h$.} For $s\ge0$, let
\begin{equation}\label{Def-ah}
 a_h=\frac1{2h}
 \left(\sum_{s=h}^{2h-1}\mathbf{1}_{S_s(o)}
       -\sum_{s=2h}^{3h-1}\mathbf{1}_{S_s(o)}\right),
\end{equation}
where $S_s(o)$ is defined as in \eqref{Def-descendantlayer}. Then $a_h$ has integral zero, is supported in $R_h^{3h}(o)$, and satisfies
\[
 \|a_h\|_\infty=\frac1{2h}
 =\frac1{\mu(R_h^{3h}(o))}.
\]
Thus $a_h$ is an atom and $\|a_h\|_{\hardy}\le1$.

\textbf{Step 2: construction of a $\mathrm{BMO}$ function.} Next, we want to construct a \(\mathrm{BMO}\) function and use a duality argument. 
Let
\[
 w_k=p^k(o),\,\,\,\gamma=\{w_k:k\ge0\},\,\,\,b(x)=\log(1+d(x,\gamma)).
\]
Firstly, we prove that
\begin{equation}\label{Ieq-bmo}
\|b\|_{\mathrm{BMO}(\mu)}\le A_q,
\end{equation}
for some $A_q>0$. By \eqref{Ieq-bmocomparison}, it suffices to find, for each admissible trapezoid $R$, a constant $c_R$ for which 
\[
\frac1{\mu(R)}\sum_{x\in R}| b(x)-c_R|\mu(x)
\]
is bounded independently of $R$.

Singletons have zero mean oscillation. Consider the admissible trapezoid
\[ 
R=R_{h'}^{h''}(y_R) = \{y\in\TT:\ y\le y_R,\ h'\le d(y,y_R)<h''\}, 
\]
and choose
\[
c_R=\log(1+d(y_R,\gamma)+h'').
\]
If $y_R\notin\gamma$, the path from any descendant $x$ of $y_R$ to $\gamma$ passes through $y_R$. Consequently,
\[
d(x,\gamma)=d(y_R,\gamma)+d(x,y_R)\ge d(y_R,\gamma)+h'.
\]
For $x\in R$, this gives
\[
0\le c_R-b(x)\le\log\frac{1+d(y_R,\gamma)+h''} {1+d(y_R,\gamma)+h'}\le\log12.
\]
Hence the required average is uniformly bounded. Now suppose that there exists some \(k\) such that $y_R=w_k\in\gamma$. In this case, $c_R=\log(1+h'')$. Fix a layer
\[
 S_t(y_R)=\{x\le y_R:d(x,y_R)=t\},\,\,\,h'\le t<h''.
\]
We count the vertices in this layer according to their distance from $\gamma$. Set
\[
 E_j(t)=\{x\in S_t(y_R):d(x,\gamma)=j\},\,\,\, 0\le j\le t.
\]
These sets partition $S_t(y_R)$. And 
\begin{equation}\label{Value-bonEj}
b(x)=\log(1+j) \,\,\,\mathrm{on}\,\, E_j(t).
\end{equation}

If $j=0$, if $t\le k$, then $E_0(t)=\{w_{k-t}\}$; if $t>k$, this set is empty. Thus there is at most one such vertex. Since $b=0$ on $\gamma$, its contribution to $|b-c_R|$ is $c_R$.

Next let $j\ge1$ and suppose that $E_j(t)$ is nonempty. For a vertex in $E_j(t)$, the path from $y_R$ first follows $\gamma$ for $t-j$ edges and then leaves it
for the remaining $j$ edges. The point where it leaves $\gamma$ is therefore $v=w_{k-t+j}$. If $v\ne o$, there are $q-1$ choices for the first edge leaving $\gamma$, followed by $q$ choices for each of the remaining $j-1$ edges. Hence $\#E_j(t)=(q-1)q^{j-1}$. If $v=o$, all $q$ children of $o$ lie outside $\gamma$, so instead $\#E_j(t)=q^j$. In particular, the bound
\begin{equation}\label{Bound-Ejcount}
 \#E_j(t)\le q^j
\end{equation}
holds for every $0\le j\le t$, including $j=0$ and the cases where $E_j(t)$ is empty.
Combining with \eqref{Value-bonEj} and \eqref{Bound-Ejcount},
\begin{equation}\label{Sum1}
\begin{aligned}
     &\frac1{\mu(R)}\sum_{x\in S_t(y_R)}|b(x)-c_R|\mu(x) = \sum_{j=0}^{t}\bigl[c_R-\log(1+j)\bigr]\frac{\mu(E_j(t))}{\mu(R)}\\
     &\le\sum_{j=0}^{t}\bigl[c_R-\log(1+j)\bigr]\frac{\#E_j(t) q^{-t}\mu(y_R)}{(h''-h')\mu(y_R)}\le \frac{1}{h''-h'}\sum_{j=0}^{t}\bigl[c_R-\log(1+j)\bigr]q^{j-t}.
\end{aligned}
\end{equation}
The term $j=0$ includes the contribution of the vertex on $\gamma$, when it exists. Observe that 
\begin{align*}
 c_R-\log(1+j)=\log\frac{1+h''}{1+t}+\log\frac{1+t}{1+j}\le \log12+\frac{t-j}{1+j}\le \log12+(t-j).
\end{align*}
Putting it back to \eqref{Sum1}, we can deduce
\begin{align*}
 \frac1{\mu(R)}\sum_{x\in S_t(y_R)}|b(x)-c_R|\mu(x) 
 &\le \frac{1}{h''-h'} \left((\log12)\sum_{j=0}^t q^{j-t}+\sum_{j=0}^{t}(t-j)q^{j-t}\right)\\
 &=\frac{1}{h''-h'} \left((\log12) \frac{q-q^t}{q-1}+\sum_{m=0}^{t}m q^{-m}\right)\\
 &\le\frac{1}{h''-h'} \left((\log12)\frac{q}{q-1}+\frac{q}{(q-1)^2}\right).
\end{align*}
This bound is uniform in both $t$ and $k$. Averaging over the layers therefore gives
\[
 \frac1{\mu(R)}\sum_{x\in R}|b(x)-c_R|\mu(x)\le\frac{q}{(q-1)^2}\left((q-1)\log 12+1\right)\eqdef \frac12 A_q.
\]
Together with the first case, this proves \eqref{Ieq-bmo}.

Finally, let $b_M=\min\{b,M\}$. So by \eqref{Ieq:bmotruncation},
\begin{equation}\label{Ieq-bMBMO}
\sup_{M>0}\|b_M\|_{\mathrm{BMO}(\mu)}\le2\|b\|_{\mathrm{BMO}(\mu)}\le 2A_q.
\end{equation}

\textbf{Step 3: duality argument.} For a fixed $h\ge 2$, define
\begin{equation}\label{Def-uhgh}
 u_h(x)\defeq \LL^{-1/2}a_h, \,\,\,g_h(x)\defeq\RR a_h(x)=\nabla u_h(x)=u_h(x)-u_h(p(x)),
\end{equation}
according to the definition of $\RR$ in \eqref{Def-verRiesz}.
We know $\RR$ is $\hardy\to L^1(\mu)$ bounded. In particular, $g_h\in L^1(\mu)$, so its pairings with any $b_M\in L^\infty(\mu)$ are well defined. 

We shall prove that, for every integer $h\ge2$, there exists $M_h>0$, which will be decided later, such that
\begin{equation}\label{Ieq-wewant}
c_q'\log h-C_q'\le\left|\sum_{x\in\TT}g_h(x)b_{M_h}(x)\mu(x)\right|\leq 2A_q\|g_h\|_{\hardy},
\end{equation}
where $c_q',C_q'>0$ are independent of $h$. If $\RR$ were bounded on $\hardy$, then
$g_h=\RR a_h$ and $\|a_h\|_{\hardy}\le1$ would imply
\[
\sup_{h\ge2}\|g_h\|_{\hardy}<\infty,
\]
which would contradict the preceding estimate as $h\to\infty$. The upper bound in \eqref{Ieq-wewant} follows from \eqref{Ieq-bMBMO}. It therefore remains to establish the lower bound for a suitable choice of $M_h$.

For $x,y\in\TT$, let $x\wedge y$ denote their nearest common ancestor. According to \eqref{Kernel-LL1}, \eqref{For-G}, we can also interpret the kernel $K_{\LL^{-1/2}}$ as 
\begin{equation}\label{Kernel-LL}
K_{\LL^{-1/2}}(x,y)=\frac{A(d(x,y))}{\mu(x\wedge y)}.
\end{equation}
where we define $A$ for all real $t\ge0$ by
\begin{equation}\label{Def-A}
 A(t)=\frac{2\sqrt2}{\pi}\sum_{j=0}^{\infty}q^{-j}\frac{t+2j+1}{(t+2j+1)^2-1/4}.
\end{equation}
For integer arguments, this is the stated kernel formula.

Put
\[
 T_o=\{x:x\le o\},\,\,\,
 D_k=\{x:x\le w_k\}\setminus\{x:x\le w_{k-1}\},
 \,\,\,k\ge1.
\]
The sets $D_k$ form a disjoint partition of $\TT\setminus T_o$. We first compute functions $u_h$ and $g_h$ on $\TT\setminus T_o$.

Fix $x\in D_k$ and write $r=d(x,w_k)$. The vertex $x$ is either $w_k$ itself or lies in a descendant branch of $w_k$ other than the one containing $o$. Hence, for every $y\in S_s(o)$, the path from $x$ to $y$ passes through $w_k$ and then $o$. In particular,
\[
 x\wedge y=w_k,\,\,\,d(x,y)=r+k+s.
\]
Thus the kernel is constant as $y$ ranges over $S_s(o)$ from \eqref{Kernel-LL}. Since $\mu(w_k)=q^k$ and $\mu(S_s(o))=1$, combining this with \eqref{Def-ah}, \eqref{Def-uhgh} and \eqref{Kernel-LL}, we obtain
\begin{align*}
 u_h(x)&=\frac{1}{2h}\sum_{s=h}^{2h-1}\sum_{y\in S_s(o)}K_{\LL^{-1/2}}(x,y)\mu(y)-\frac{1}{2h}\sum_{s=2h}^{3h-1}\sum_{y\in S_s(o)}K_{\LL^{-1/2}}(x,y)\mu(y)\\
 &=\frac{1}{2h}\sum_{s=h}^{2h-1}\sum_{y\in S_s(o)}\mu(y)\frac{A(k+r+s)}{\mu(w_k)}-\frac{1}{2h}\sum_{s=2h}^{3h-1}\sum_{y\in S_s(o)}\mu(y)\frac{A(k+r+s)}{\mu(w_k)}\\
 &=\frac{q^{-k}}{2h}\sum_{s=h}^{2h-1}\bigl[A(k+r+s)-A(k+r+s+h)\bigr].
\end{align*}
Define
\begin{equation}\label{Def-F}
 F_h(t)=\frac1{2h}\sum_{s=h}^{2h-1}\bigl[A(t+s)-A(t+s+h)\bigr], \,\,\,t\ge0.
\end{equation}
With this notation,
\[
 u_h(x)=q^{-k}F_h(k+r).
\]
If $r\ge1$, then $p(x)\in D_k$ and $d(p(x),w_k)=r-1$. Applying the same formula
at $p(x)$, and recalling the definition of $g_h$ in \eqref{Def-uhgh}, for $x\in D_k$ with $r\ge1$,
\begin{equation}\label{Idn-gh}
 \begin{aligned}
g_h(x)=q^{-k}\bigl[F_h(k+r)-F_h(k+r-1)\bigr].
 \end{aligned}
\end{equation}

We next establish the properties of $F_h$ needed to estimate this difference.
The series in \eqref{Def-A}, together with its first two derivative series, converges uniformly on $[0,\infty)$: each summand in each series is bounded
by a constant times $q^{-j}$. Termwise differentiation is therefore justified.
Every summand in \eqref{Def-A} is positive, strictly decreasing, and strictly convex.
Moreover, keeping the $j=0$ term gives lower bounds, while summing the geometric series gives upper bounds. It follows that
\begin{equation}\label{Est-AAderi}
 A(t)\simeq_q\frac1{1+t},\,\,\,-A'(t)\simeq_q\frac1{(1+t)^2},\,\,\, A''(t)>0.
\end{equation}
Thus,
\begin{equation}\label{Est-F}
 F_h(t)=\frac1{2h}\sum_{s=h}^{2h-1}\int_0^h[-A'(t+s+v)]\dd v>0.
\end{equation}
Differentiating this identity gives
\[
 F_h'(t)=-\frac1{2h}\sum_{s=h}^{2h-1}\int_0^h A''(t+s+v)\dd v<0.
\]
Thus $F_h$ is positive and strictly decreasing. Finally, for $h\le s<2h$ and $0\le v\le h$, we have $1+t+s+v\simeq t+h$. Hence by \eqref{Est-AAderi},
\[
 \int_0^h[-A'(t+s+v)]\dd v\simeq_q\frac{h}{(t+h)^2}.
\]
Combining this with \eqref{Est-F}, we conclude
\begin{equation}\label{Est-Fbounds}
 F_h(t)\simeq_q\frac{h}{(t+h)^2},\,\,\, t\ge0,\,\, h\ge2.
\end{equation}

For $k\ge1$ and $r\ge1$, set
\[
E_{k,r}\defeq\{x\in D_k:d(x,w_k)=r\}.
\]
By definition, 
\begin{equation}\label{Count-Ekr}
\#E_{k,r}=(q-1)q^{r-1}.
\end{equation}
Since $\ell(w_k)=k$, every $x\in E_{k,r}$ satisfies $\ell(x)=k-r$ and hence $\mu(x)=q^{k-r}$. Moreover, for $x\in E_{k,r}$, $d(x,\gamma)=r$. It follows that
\begin{equation}\label{Idn-bM}
b_M(x)=\min\{\log(1+r),M\}\eqdef\phi_M(r).
\end{equation}
The remaining vertex $w_k$ corresponds to $r=0$ and contributes nothing to the pairing, since $b_M(w_k)=0$. Thus, putting \eqref{Idn-gh}, \eqref{Count-Ekr} and \eqref{Idn-bM} together, we can deduce 
\begin{equation*}
\begin{aligned}
 I_{\mathrm{out}}(h,M) &\defeq\sum_{x\in\TT\setminus T_o}g_h(x)b_{M}(x)\mu(x)=\sum_{k\ge1}\sum_{r\ge1}\sum_{x\in E_{k,r}}g_h(x)b_{M}(x)\mu(x)\\
 &=\sum_{k\ge1}\sum_{r\ge1} q^{-k}[F_h(k+r)-F_h(k+r-1)]\phi_M(r) q^{k-r} \#E_{k,r}\\
 &=(1-q^{-1})\sum_{k\ge1}\sum_{r\ge1}[F_h(k+r)-F_h(k+r-1)]\phi_M(r)\\
 \end{aligned}
\end{equation*}
By telescopic sum in $k$ for fixed $r$, we have 
 \begin{equation}\label{Iout}
 I_{\mathrm{out}}(h,M)=-(1-q^{-1})\sum_{n\ge1}F_h(n)\phi_M(n).
\end{equation}
Choose
\begin{equation}\label{Def-Mh}
 M_h=\log(2h+1).
\end{equation}
For $h\le n\le2h$, we have $\phi_{M_h}(n)=\log(n+1)$. Keeping only these terms in
\eqref{Iout} and using \eqref{Est-Fbounds}, by the nonnegativity of $F_h(n)\phi_{M_h}(n)$, we obtain
\begin{equation}\label{Ieq-outlower}
 I_{\mathrm{out}}(h,M_h)\le-c_q'\sum_{n=h}^{2h}\frac{h\log(n+1)}{(h+n)^2} \le-c_q'\log h,
\end{equation}
for some positive $c_q'$ independent of $M$ and $h$.

It remains to estimate 
\[
I_{\mathrm{in}}(h,M_h)\defeq \sum_{x\in T_o}g_h(x) b_{M_h}(x)\mu(x).
\]
We will use the $L^2$ boundedness of $\RR$ near the support of the atom $a_h$ and kernel estimates farther away. Set
\[
\widetilde b_{h,M_h}(x)=b_{M_h}(x)-\phi_{M_h}(h),\,\,\,E_h=\bigcup_{t=0}^{6h}S_t(o).
\]
We decompose
\begin{equation}\label{Idn-indecom}
\begin{aligned}
&I_{\mathrm{in}}(h,M_h)\\
&=\sum_{x\in E_h}g_h(x)\widetilde b_{h,M_h}(x)\mu(x)
 +\sum_{x\in T_o\setminus E_h}g_h(x)\widetilde b_{h,M_h}(x)\mu(x)
 +\phi_{M_h}(h)\sum_{x\in T_o}g_h(x)\mu(x)\\
 &\eqdef I_{\mathrm{in}}^{1,1}+I_{\mathrm{in}}^{1,2}+I_{\mathrm{in}}^{2}.
\end{aligned}
\end{equation}

We first estimate $I_{\mathrm{in}}^{1,1}$. The $L^2$ boundedness of $\RR$ gives
\[
\|g_h\|_{L^2(\mu)}\lesssim\|a_h\|_{L^2(\mu)}\leq\|a_h\|^{1/2}_{L^1(\mu)}\|a_h\|^{1/2}_{L^\infty(\mu)}=\frac1{\sqrt{2h}}.
\]
For $x\in S_t(o)$, by \eqref{Idn-bM} and \eqref{Def-Mh}, we have 
\begin{equation}\label{Est-tildebMh}
\begin{aligned}
    \left|\widetilde b_{h,M_h}(x)\right|&=\left|\phi_{M_h}(t)-\phi_{M_h}(h)\right|=\left|\log\frac{1+t}{1+h}\right|\cdot\mathbf{1}_{\{0\le t< 2h\}}+\log \frac{1+2h}{1+h}\cdot\mathbf{1}_{\{t\ge2h\}}.
\end{aligned}
\end{equation}
We obtain
\begin{equation}\label{Importantestimate}
\begin{aligned}
\|\widetilde b_{h,M_h}\mathbf1_{E_h}\|_{L^2(\mu)}^2
\le &\sum_{t=0}^{h-1}\left(\log\frac{1+h}{1+t}\right)^2+\sum_{t=h}^{6h}\left(\log\frac{1+t}{1+h}\right)^2
\lesssim \sum_{j=1}^{h}\log^2 \frac{1+h}{j}+h.
\end{aligned}
\end{equation}
Since the function $x\mapsto \left(\log \frac{1+h}{x}\right)^2$ is nonnegative and decreasing, then we deduce that 
\[
\sum_{j=1}^{h}\log^2 \frac{1+h}{j}\le \int_0^h\log^2 \frac{1+h}{x}\dd x\le\int_0^{h+1}\log^2 \frac{1+h}{x}\dd x.
\]
Using change of variable $x=(1+h)\ee^{-r}$, then 
\[
\int_0^{h+1}\log^2 \frac{1+h}{x}\dd x=(1+h)\int_0^\infty r^2\ee^{-r}\dd r=2(1+h)\lesssim h.
\]
Plugging it into \eqref{Importantestimate}, it follows that 
\[
\|\widetilde b_{h,M_h}\mathbf1_{E_h}\|_{L^2(\mu)}\lesssim h^{1/2}.
\]
Hence the Cauchy--Schwarz inequality yields
\begin{equation}\label{Ieq-innear}
I_{\mathrm{in}}^{1,1}\le\|g_h\|_{L^2(\mu)}\|\widetilde b_{h,M_h}\mathbf1_{E_h}\|_{L^2(\mu)}\lesssim_q 1.
\end{equation}

We next estimate $I_{\mathrm{in}}^{1,2}$. Fix $x\in S_t(o)$ with $t>6h$, then by \eqref{Kernel-LL} and \eqref{Def-uhgh}, we can conclude that 
\begin{equation}\label{Est-I12gh}
    \begin{aligned}
        g_h(x)=&\sum_{y\in \bigcup_{h\le s<3h}S_s(o)}\left[K_{\LL^{-1/2}}(x,y)-K_{\LL^{-1/2}}(p(x),y)\right]a_h(y)\mu(y)\\
        =&\sum_{y\in \bigcup_{h\le s<3h}S_s(o)}\left[A(d(x,y))-A(d(x,y)-1)\right]a_h(y)\frac{\mu(y)}{\mu(x\wedge y)},
    \end{aligned}
\end{equation}
since for every $y\in S_s(o)$ with $h\le s<3h$,
\[
p(x)\wedge y=x\wedge y,\,\,\,d(p(x),y)=d(x,y)-1.
\]
Moreover, for this kind of \(y\),
\[
d(x,y)\ge t-s>\frac t2.
\]
The estimates for $A$ in \eqref{Est-AAderi} therefore imply
\begin{equation}\label{Est-AminusA}
|A(d(x,y))-A(d(x,y)-1)|\lesssim_q\frac1{(1+t)^2}.
\end{equation}
For $s<t$, we also have the elementary bound
\begin{equation}\label{Ieq-ancestorcount}
\sum_{y\in S_s(o)}\frac{\mu(y)}{\mu(x\wedge y)}\le s+1.
\end{equation}
Combining with \eqref{Est-I12gh}, \eqref{Est-AminusA} and \eqref{Ieq-ancestorcount}, 
\begin{equation}\label{Est-I12ghab}
\begin{aligned}
|g_h(x)|&\lesssim_q\frac1{2h(1+t)^2}\sum_{s=h}^{3h-1}(s+1)\lesssim_q\frac{h}{(1+t)^2}.
\end{aligned}
\end{equation}
Consequently, combining with \eqref{Est-tildebMh} for $t> 6h$, 
\begin{equation}\label{Ieq-infar}
I_{\mathrm{in}}^{1,2}\lesssim_q\sum_{t>6h}\frac{h}{(1+t)^2}\log\frac{1+2h}{1+h}\lesssim_q 1.
\end{equation}

It remains to control $I_{\mathrm{in}}^2$. By symmetry, $u_h$ is constant on each $S_t(o)$; denote this value by $U_h(t)$. Similar with \eqref{Est-I12ghab}, we can obtain
\begin{equation}\label{Est-uh}
    |U_h(t)|=|u_h(x)|\lesssim_q\frac{1}{2h(1+t)}\sum_{s=h}^{3h-1}(s+1)\lesssim_q\frac{h}{1+t}.
\end{equation}
For every integer $N\ge1$, telescoping gives
\begin{align*}
\sum_{t=0}^{N}\sum_{x\in S_t(o)}g_h(x)\mu(x)
&=U_h(0)-u_h(p(o))+\sum_{t=1}^{N}[U_h(t)-U_h(t-1)]\\&=U_h(N)-u_h(p(o)).
\end{align*}
We already know that $g_h\in L^1(\mu)$, and \eqref{Est-uh} gives $U_h(N)\to0$ as $N\to\infty$. Therefore,
\[
\sum_{x\in T_o}g_h(x)\mu(x)=-u_h(p(o))=-q^{-1}F_h(1).
\]
and from \eqref{Est-Fbounds}, we conclude that
\begin{equation}\label{Ieq-inconstant}
\phi_{M_h}(h)\left|\sum_{x\in T_o}g_h(x)\mu(x)\right|\lesssim_q\frac{\log(1+h)}h\lesssim_q1.
\end{equation}

Combining \eqref{Ieq-innear}, \eqref{Ieq-infar} and \eqref{Ieq-inconstant}, we obtain
\begin{equation}\label{Iin}
|I_{\mathrm{in}}(h,M_h)|\le C_q',\,\,\, h\ge2,
\end{equation}
with some $C_q'$ independent of $h$ and $M$.

Combining \eqref{Ieq-outlower} and \eqref{Iin}, we deduce the lower bound in \eqref{Ieq-wewant} with setting $C_q=\frac{C_q'}{2A_q}$ and $c_q=\frac{c_q'}{2A_q}$.
\end{proof}

\section{Failure of weak type \texorpdfstring{$(1,1)$}{(1,1)} for the adjoint vertical Riesz transform}\label{Sec-adjver}

Having established that $\RR$ is not bounded on $\hardy$, we now turn to its adjoint. We show that $\RR^*$ is not of weak type $(1,1)$, providing a discrete counterpart to Hormozi's result for the adjoint vertical Riesz transform on $ax+b$ groups; see \cite[Theorem~1.1]{H26}.

\begin{proof}[Proof of Theorem~\ref{Thm-mainadjver}]
Since $\RR$ is of weak type $(1,1)$ by \cite[Theorem~2.3]{HS03} and \cite{ATV20}, it suffices to prove that $\RR-\RR^*$ is not of weak type $(1,1)$. By \cite[Eqs.~(3.4) and (4.4)]{LMSTV23}, for every finitely supported
function $f$,
\begin{equation}\label{Idn-adjskew}
(\RR-\RR^*) f(x) =\sum_{n=1}^{\infty}\widetilde k_{\ZZ}(n)\left((\Sigma^*)^n f(x)-\Sigma^n f(x)\right),
\end{equation}
here we used the property that $\widetilde k_{\ZZ}$ is odd by \eqref{For-G}.
And the notation $(\Sigma^*)^n f(x)$ means the average of $f$ among the $n$-generation descendants of vertex $x$, while $\Sigma^n f(x)=f(p^n(x))$.

Let $N$ be a sufficiently large positive integer. Define
\begin{equation}\label{Choice-m}
m=\lfloor \log_q(2N)\rfloor+1,\,\,\,\xi_N\defeq p^{2N}(o).
\end{equation}
For $0\le h\le N$ and every vertex of $S_h(\xi_N)$, choose one $m$-generation descendant, then we denote this set by $\Gamma_h$. Thus, 
\begin{equation}\label{Idn-Gammameasure}
\mu(\Gamma_h)=q^{-m}\sum_{x\in S_{h}(\xi_N)}\mu(x)=q^{-m}\mu(\xi_N)=q^{-m+2N}.
\end{equation}
Define
\begin{equation}\label{Def-fN}
f_N=\frac{q^{m-2N}}{(N+1)}\sum_{h=0}^{N}\mathbf1_{\Gamma_h},
\end{equation}
so that
\[
\|f_N\|_{L^1(\mu)}=1.
\]
It suffices to find sets $E_N\subseteq\TT$ and numbers $\lambda_N>0$ such that
\[
(\RR-\RR^*)f_N(x)\ge \lambda_N, \,\,\,\forall x\in E_N,
\]
and 
\[
\lambda_N\mu(E_N)\to \infty, \,\,\,\mathrm{as}\,\,N\to\infty.
\]
Indeed, 
\begin{equation}\label{Weaktype11}
\begin{aligned}
\|(\RR-\RR^*)f_N\|_{L^{1,\infty}(\mu)}\ge\lambda_N\,
\mu\bigl(\{x\in\TT:|(\RR-\RR^*) f_N(x)|>\lambda_N\}\bigr)\ge \lambda_N\mu(E_N)
\to\infty.
\end{aligned}
\end{equation}
We first estimate the contribution of the descendant averages. For fixed $x\in S_h(\xi_N)$ for $0\le h\le N$, every set $S_n(x)$ only contains $q^{n-m}$ selected vertices for $m\le n\le N+m-h$. 
\begin{equation}\label{Sigmaadn}
    (\Sigma^*)^nf_N(x)=q^{-n}q^{n-m}\frac{q^{m-2N}}{(N+1)}=\frac{q^{-2N}}{N+1}, \,\,\,m\le n\le N+m-h.
\end{equation}
For $n>N+m-h$, it contains $0$ selected vertices. Since $f_N\ge0$ and $\widetilde k_{\ZZ}(n)>0$ for $n\ge1$, all terms are nonnegative. Hence
\begin{equation}\label{Sum-tildekn}
\sum_{n=1}^{\infty}\widetilde k_{\ZZ}(n)(\Sigma^*)^n f_N(x)
\ge\frac{q^{-2N}}{N+1}\sum_{n=m}^{N+m-h}\widetilde k_{\ZZ}(n).
\end{equation}

We now consider the ancestor terms. Indeed, for any $x\in\TT$, 
\[
\Sigma^n f_N(x)\ne 0\Leftrightarrow f_N(p^n(x))\neq 0\Leftrightarrow p^n(x)\in \bigcup_{0\leq h\leq N} \Gamma_h.
\] 
To make these terms vanish, we remove from each layer all vertices having a selected vertex as a strict ancestor. 
More precisely, for $0\le h\le N$, define
\[
F_h\defeq S_h(\xi_N)\setminus\bigcup_{j=0}^{N}\bigcup_{y\in\Gamma_j}
\{z\in\TT:z\le y,\ z\ne y\}.
\]
For every $x\in F_h$, none of its strict ancestors belongs to the support of $f_N$. Consequently,
\begin{equation}\label{Sigman}
\Sigma^n f_N(x)=f_N(p^n(x))=0,\,\,\, x\in F_h,\,\,n\ge1.
\end{equation}
Moreover, recall that $\Gamma_j$ consists of $q^j$ selected vertices in $S_{j+m}(\xi_N)$. Only indices $j$ satisfying $j+m<h$ contribute to the removal of vertices from $S_h(\xi_N)$.
Fix such an index $j$. Each vertex of $\Gamma_j$ has exactly $q^{h-j-m}$ descendants in $S_h(\xi_N)$, and every descendant has measure $q^{-h}\mu(\xi_N)$. Thus,
\[
\mu\left(\bigcup_{y\in\Gamma_j}
\{z\in S_h(\xi_N):z\le y,\ z\ne y\}\right)\leq \#\Gamma_j\cdot q^{h-j-m}\cdot q^{-h}\mu(\xi_N)=q^{-m+2N}.
\]
There are at most $N$ indices satisfying $0\le j<h-m$. Although the sets removed by different $\Gamma_j$ may overlap, their union therefore satisfies
\[
\mu\left(\bigcup_{j=0}^{N}\bigcup_{y\in\Gamma_j}
\{z\in S_h(\xi_N):z\le y,\ z\ne y\}\right)\le Nq^{2N-m}.
\]
Consequently, our choice of $m$ in \eqref{Choice-m} gives
\begin{equation}\label{measure-Fh}
\mu(F_h)\ge (1-Nq^{-m})q^{2N}>\frac12q^{2N},\,\,\,0\le h\le N.
\end{equation}

On every $F_h$ for $0\le h\le N$, \eqref{Sigmaadn} and \eqref{Sigman} hold, yielding \eqref{Sum-tildekn}. Combining with \eqref{Idn-adjskew}, for every $x\in F_h$,
\[
(\RR-\RR^*) f_N(x)=\sum_{n=1}^{\infty}\widetilde k_{\ZZ}(n)(\Sigma^*)^n f_N(x)
\ge\frac{q^{-2N}}{N+1}\sum_{n=m}^{N+m-h}\widetilde k_{\ZZ}(n).
\]

Set
\[
E_N=\bigcup_{h=0}^{\lfloor N/2\rfloor}F_h,\,\,\,L_N=\sum_{n=m}^{\lfloor N/2\rfloor}
\widetilde k_{\ZZ}(n),\,\,\,\lambda_N=\frac{q^{-2N}}{N+2}L_N.
\]
For $0\le h\le\lfloor N/2\rfloor$, we have $N+m-h\ge\lfloor N/2\rfloor$, and hence
\[
(\RR-\RR^*)f_N(x)\ge\frac{q^{-2N}}{N+1}L_N>\lambda_N,\,\,\, x\in E_N.
\]
Additionally, \eqref{measure-Fh} yields
\[
\mu(E_N)=\sum_{h=0}^{\lfloor N/2\rfloor}\mu(F_h)\ge\frac{\lfloor N/2\rfloor+1}{2}\,q^{2N}\ge\frac{N+1}{4}\,q^{2N}.
\]
Consequently,
\[
\lambda_N\mu(E_N)\ge\frac14\frac{N+1}{N+2}L_N>\frac18L_N.
\]
Finally, since $\widetilde k_{\ZZ}(n)\gtrsim n^{-1}$ by \eqref{For-G}, and $m=O_q(\log N)$,
\[
L_N\gtrsim\sum_{n=m}^{\lfloor N/2\rfloor}\frac1n\ge\log\frac{\lfloor N/2\rfloor+1}{m}
\to\infty.
\]
Thus \eqref{Weaktype11} holds.
\end{proof}

\section{Horizontal Riesz transforms on \texorpdfstring{$\hardy$}{H1}} \label{Sec-hor}

We next turn to horizontal Riesz transforms, whose behaviour on $\hardy$ depends on the choice of discrete derivative. We first show that the scalar construction of \cite[Section~5]{LMSTV23} fails to preserve cancellation on atoms. We then introduce a horizontal derivative with cancellation on each sibling set. The associated Riesz transform retains all the boundedness of the scalar horizontal transform. Moreover, we prove that it is also bounded on $\hardy$ itself.
\subsection{The scalar horizontal Riesz transform}
In this subsection, we will prove that $\RR_\var$ is $\hardy\to L^1(\mu)$ bounded, but one can not strengthen this boundedness to $\hardy$ to itself.
\begin{proof}[Proof of Theorem~\ref{Thm-varH1}]
According to \cite[Theorem 1.2]{S22}, the discrete Riesz transform \(\RR=\nabla \LL^{-1/2}\) is bounded from \(\hardy\) to \(L^1(\mu)\). In order to show the \(\hardy\) to \(L^1(\mu)\) boundedness of \(\RR_{\var}\), using \eqref{Rel-derscalarvert}, it suffices to show that \(\nabla_\var\) is bounded from $L^1(\mu)$ to $L^1(\mu)$, which is obvious from \eqref{Est-nablavar}. 

However, the \(L^1(\mu)\)-estimate alone is not enough to obtain \(\hardy\)-boundedness. The additional point is that this operator lacks the required cancellation. Next we construct a \((1,\infty)\) atom \(a\) to show 
\[
\sum_{x\in\TT} \RR_\var a(x)\mu(x)\neq 0.
\]

Since \(\var\) is nonzero, choose $\xi\in\TT$ such that $\var(\xi)\neq0$ and set $y_R=p(\xi)$. The cancellation condition
\[
\sum_{z\in s(y_R)}\var(z)=0
\]
ensures that there exists $\eta\in s(y_R)$ with $\var(\eta)\neq\var(\xi)$; otherwise, the sum would equal $q\var(\xi)\neq0$.

Consider the admissible trapezoid
\[
R=\{y\le y_R:\ 1\le d(y,y_R)<2\}=s(y_R).
\]
Define
\[
a=\frac1q\left(\frac{\delta_{\{\xi\}}}{\mu(\xi)}-\frac{\delta_{\{\eta\}}}{\mu(\eta)}\right).
\]
Since \(\mu(\xi)=\mu(\eta)\), we have
\[
\sum_{y\in R}a(y)\mu(y)=0,\,\,\,\|a\|_\infty=(\mu(R))^{-1}.
\]
Thus \(a\) is a \((1,\infty)\)-atom supported in \(R\).

We now compute the cancellation of \(\RR_\var a\). From formulas \eqref{Kernel-LL1} and \eqref{For-G},
\[
\begin{aligned}
\RR_\var a(x)
&=\frac1q\sum_{z\in s(x)}\var(z)\LL^{-1/2}a(z)
=\frac1{q^2}\sum_{z\in s(x)}\var(z)\left(K_{\LL^{-1/2}}(z,\xi)-K_{\LL^{-1/2}}(z,\eta)\right)\\
&=\frac1{q^2}q^{-(\ell(x)-1)/2}q^{-\ell(\xi)/2}\sum_{z\in s(x)}\var(z)(G(d(z,\xi))-G(d(z,\eta))).
\end{aligned}
\]
Hence,
\begin{equation}\label{For-sum}
\begin{aligned}
\sum_{x\in\TT} \RR_\var a(x)\mu(x)
=&\frac1{q}\sum_{x\in\TT}q^{(\ell(x)-\ell(\xi)-1)/2} \sum_{z\in s(x)}\var(z) (G(d(z,\xi))-G(d(z,\eta))).
\end{aligned}
\end{equation}

We claim that all terms with $x\neq y_R$ vanish. First suppose that $x\neq y_R$ belongs to neither of the descendant subtrees rooted at $\xi$ and $\eta$. Then for every $z\in s(x)$, $d(z,\xi)=d(z,\eta)$, so
\[
G(d(z,\xi))-G(d(z,\eta))=0.
\]
Next suppose that $x\le\xi$ and put $r=d(x,\xi)$. Every $z\in s(x)$ satisfies
\[
d(z,\xi)=r+1,\,\,\, d(z,\eta)=r+3.
\]
Consequently,
\[
\sum_{z\in s(x)}\var(z)(G(d(z,\xi))-G(d(z,\eta)))=(G(r+1)-G(r+3))\sum_{z\in s(x)}\var(z)=0.
\]
The case $x\le\eta$ is analogous.

It remains to consider $x=y_R$. In $s(x)$, only $z=\xi$ and $z=\eta$ contribute, and hence
\[
\sum_{z\in s(y_R)}\var(z)(G(d(z,\xi))-G(d(z,\eta)))=(\var(\xi)-\var(\eta))(G(0)-G(2)).
\]
Using the formulas (\ref{For-Gnn+2}) and \eqref{For-sum}, we obtain
\[
\sum_{x\in\TT}\RR_\var a(x)\mu(x)
=\frac1q\bigl(\var(\xi)-\var(\eta)\bigr)\bigl(G(0)-G(2)\bigr)
=\frac1q\,(\var(\xi)-\var(\eta))\widetilde{k}_{\ZZ}(1)\neq 0.
\]

Thus \(\RR_\var a\) does not have zero integral. Since every element of
\(\hardy\) has cancellation, this shows that \(\RR_\var a\notin \hardy\). Therefore \(\RR_\var\) cannot be bounded from \(\hardy\) to itself.
\end{proof}

Motivated by this failure of cancellation, we introduce a horizontal derivative based on deviations from sibling averages.
\subsection{A horizontal Riesz transform with sibling cancellation}
Recall that the horizontal derivative and its associated Riesz transform are defined by
\begin{equation}\label{Def-newhorRiesz}
\nabla_h f(x)=\frac1{q}\left(f(x)-\frac1{q-1}\sum_{z\in b(x)}f(z)
\right),\,\,\,\RR_h=\nabla_h\LL^{-1/2},
\end{equation}
which has the cancellation across the \(q\) children of each vertex. 
\begin{equation}\label{Idn-cancellationq}
\sum_{x\in s(y)}\RR_h f(x)\mu(x)=0.
\end{equation}

To establish the boundedness properties of $\RR_h$, we next express it in terms of finitely many scalar horizontal Riesz transforms $\RR_{\var_j}$ with suitably chosen coefficients $\var_j$. This representation allows us to
transfer the known estimates for the scalar transforms to $\RR_h$.

For each vertex \(y\), the children \(s(y)\) are ordered from left to
right:
\[
s(y)=\{y_1,\ldots,y_q\}.
\]
For \(j=1,\ldots,q\), define \(\var_j\) on each sibling set \(s(y)=\{y_1,\ldots,y_q\}\) by
\[
\var_j(y_j)=1,\,\,\,\var_j(y_i)=-\frac1{q-1},\,\,\, i\neq j,
\]
so that $\sum_{i=1}^q\var_j(y_i)=0$. Then for \(f:\TT\to\mathbb C\), by \eqref{Def-horizontalscalar}, 
\begin{equation}\label{Def-nablavarj}
\nabla_{\var_j}f(y)=\frac1q\sum_{i=1}^q\var_j(y_i)f(y_i)=\frac{1}{q-1}\left(f(y_j)-\frac1q\sum_{i=1}^q f(y_i)\right).
\end{equation}

Let
\begin{equation}\label{Def-Sj}
(M_{\var_j}g)(x)=\var_j(x)g(x),\,\,\,S_j\defeq M_{\var_j}\Sigma,
\end{equation}
where $\Sigma$ is the shift operator defined by \eqref{Def-shift}.
Now we can establish the following lemma to clarify the relation between $\RR_h$ and $\RR_{\var_j}$.
\begin{lemma}
For the $S_j$ defined above, one has
 \begin{equation}\label{Relation-RRhRRvarj}
\RR_h f=\frac{q-1}{q}\sum_{j=1}^q S_j\RR_{\var_j}f.
\end{equation}
\end{lemma}
\begin{proof}
Let \(x=y_k\in s(y)\). Because of \eqref{Def-nablavarj}, we have
\[
\begin{aligned}
\sum_{j=1}^q S_j\nabla_{\var_j}f(y_k)
=\sum_{j=1}^q\var_j(y_k)\nabla_{\var_j}f(y)
=\frac1{q-1}\sum_{j=1}^q\var_j(y_k)\left(f(y_j)-\frac1q\sum_{i=1}^q f(y_i)\right).
\end{aligned}
\]
Since
\[
\sum_{j=1}^q \varepsilon_j(y_k)=0,
\]
the average term vanishes. Hence
\[
\sum_{j=1}^q S_j\nabla_{\var_j}f(y_k)
=\frac1{q-1}\sum_{j=1}^q\var_j(y_k)f(y_j).
\]
By the definition of \(\var_j\) and \eqref{Def-newhorRiesz},
\[
\sum_{j=1}^q\var_j(y_k)f(y_j)=f(y_k)-\frac1{q-1}\sum_{\substack{j=1\\ j\ne k}}^qf(y_j)=q\nabla_h f(y_k).
\]
Consequently, we deduce
\[
\nabla_hf(y_k)=\frac{q-1}{q}\sum_{j=1}^q S_j\nabla_{\var_j}f(y_k),
\]
and hence the desired result.
\end{proof}

We now derive the boundedness properties of \(\RR_h\) from the corresponding properties of the directional Riesz transforms \(\RR_{\var_j}\), including its $L^p(\mu)$-boundedness, weak type $(1,1)$-boundedness and \(\hardy\to L^1(\mu)\) boundedness.

\begin{theorem}\label{Thm-Lpweakttype11q}
    The operators \(\RR_h\) and \(\RR_h^*\) are of weak type (1,1) and bounded on \(L^{p}(\mu)\) for every \(p\in(1,\infty)\). Moreover, \(\RR_h\) is bounded from \(\hardy\) to \(L^1(\mu)\).
\end{theorem}
\begin{proof}
Since \(|\var_j|\le1\), \(M_{\var_j}\) is a contraction on \(L^p(\mu)\) and on \(L^{1,\infty}(\mu)\). Hence, for every \(j=1,\ldots,q\), from the definition of \(S_j\) in \eqref{Def-Sj} and \eqref{shift}, the operator \(S_j\) is bounded on \(L^p(\mu)\) and on \(L^{1,\infty}(\mu)\). More precisely,
\begin{equation}\label{ID-MCLpq}
\|S_j f\|_{L^p(\mu)}=\|M_{\var_j}\Sigma f\|_{L^p(\mu)}\le\|\Sigma f\|_{L^p(\mu)}=\|f\|_{L^p(\mu)},\,\,\,\mathrm{for}\,\, 1\le p\le \infty,
\end{equation}
and
\begin{equation*}
\|S_j f\|_{L^{1,\infty}(\mu)}=\|M_{\var_j}\Sigma f\|_{L^{1,\infty}(\mu)}\leq\|\Sigma f\|_{L^{1,\infty}(\mu)}
=\|f\|_{L^{1,\infty}(\mu)} .
\end{equation*}

On the other hand, for every \(j=1,\ldots,q\), \(\RR_{\var_j}\) is bounded on \(L^p(\mu)\) for \(1<p<\infty\), and is of weak type \((1,1)\), by \cite[Section 5, p.~27]{LMSTV23}. Combining these facts with the decomposition \eqref{Relation-RRhRRvarj}, we obtain the \(L^p-\)boundedness and weak type \((1,1)\)-boundedness of \(\RR_h\). Indeed, for \(1<p<\infty\),
\[
\|\RR_h f\|_{L^p(\mu)}
\le \frac{q-1}{q}\sum_{j=1}^q\|S_j \RR_{\var_j}f\|_{L^p(\mu)} 
\le\frac{q-1}{q}\sum_{j=1}^q\|\RR_{\var_j}f\|_{L^p(\mu)}
\lesssim\|f\|_{L^p(\mu)},
\]
and 
\[
\|\RR_h f\|_{L^{1,\infty}(\mu)}
\lesssim \sum_{j=1}^q\|S_j \RR_{\var_j}f\|_{L^{1,\infty}(\mu)} 
\le\frac{q-1}{q}\sum_{j=1}^q\|\RR_{\var_j}f\|_{L^{1,\infty}(\mu)}
\lesssim\|f\|_{L^1(\mu)} .
\]

Finally, by Theorem \ref{Thm-varH1}, each \(\RR_{\var_j}\) is bounded from
\(\hardy\) to \(L^1(\mu)\). Therefore, using again \eqref{Relation-RRhRRvarj} and \eqref{ID-MCLpq} with \(p=1\), we obtain
\[
\begin{aligned}
\|\RR_h f\|_{L^1(\mu)}
\le\frac{q-1}{q}\sum_{j=1}^q\|S_j \RR_{\var_j}f\|_{L^1(\mu)}
\le\frac{q-1}{q}\sum_{j=1}^q\|\RR_{\var_j}f\|_{L^1(\mu)}
\lesssim\|f\|_{\hardy} .
\end{aligned}
\]
Consequently,
\(
\RR_h\) is \(\hardy\to L^1(\mu)
\)
bounded.

For the adjoint operator, the \(L^p\)-boundedness follows by duality from the
\(L^{p'}\)-boundedness of \(\RR_h\), where \(p'\) is the conjugate exponent of \(p\). 

It remains to prove the weak type \((1,1)\) estimate for \(\RR_h^*\). Taking adjoints with respect to the \(L^2(\mu)\)-pairing for \(\RR_h\), we obtain
\[
\RR_h^*=\frac{q-1}{q}\sum_{j=1}^q\RR_{\varepsilon_j}^*S_j^*
=\frac{q-1}{q}\sum_{j=1}^q\RR_{\varepsilon_j}^*\Sigma^*M_{\varepsilon_j}^*.
\]
It follows from \eqref{shift} that
\[
\|S_j^*f\|_{L^1(\mu)}=\|\Sigma^*M_{\var_j}^* f\|_{L^1(\mu)}
\le\|M_{\var_j}^* f\|_{L^1(\mu)}\le\|f\|_{L^1(\mu)} .
\]
Using the weak type \((1,1)\) estimate for \(\RR_{\var_j}^*\), see
\cite[Theorem 5.5]{LMSTV23}, we get
\[
\begin{aligned}
\|\RR_h^* f\|_{L^{1,\infty}(\mu)}
\lesssim \sum_{j=1}^q\|\RR_{\var_j}^*S_j^*f\|_{L^{1,\infty}(\mu)} 
\lesssim\sum_{j=1}^q\|S_j^*f\|_{L^1(\mu)}
\le\sum_{j=1}^q\|f\|_{L^1(\mu)}
=q\|f\|_{L^1(\mu)} .
\end{aligned}
\]
This proves that \(\RR_h^*\) is of weak type \((1,1)\). 
\end{proof}

Moreover, we can strengthen the above $\hardy\to L^1(\mu)$ boundedness to the boundedness on \(\hardy\), i.e., Theorem~\ref{Thm-mainhor}. We first establish a lemma.

\begin{lemma}\label{Lem-H1L1-cancellation}
    For any \(f\in L^1(\mu)\) such that \(\Sigma ^*f=0\), i.e., there is cancellation among siblings, then \(f\in \hardy\) and 
    \begin{equation}\label{Ieq-H1L1inlem}
    \|f\|_{\hardy}= \|f\|_{L^1(\mu)}.
    \end{equation}
    In other words, the kernel of \(\Sigma ^*\) in \(L^1(\mu)\) is contained in \(\hardy\) and the \(L^1(\mu)\)-norm and \(\hardy\)-norm coincide in this subspace.
\end{lemma}
\begin{proof}
By the definition of \(\Sigma^*\) in~\eqref{Def-shift}, for every \(y\in\TT\), 
\begin{equation}\label{Id-sigmaadf0}
    \Sigma^*f(y)=\frac1q\sum_{x\in s(y)}f(x)=0.
\end{equation}
For \(y\in \TT\), write \(P_{y}=s(p(y))\). Then \(P_{y}\) is an admissible trapezoid of height \(1\). By~\eqref{Id-sigmaadf0}, we always have the following decomposition:
\[
f=\sum_{y\in\TT}f(y)\mu(y)a_y,\,\,\,a_y=\frac{1}{\mu(y)}\delta_{y}-\frac1{\mu(P_y)}\delta_{P_y}.
\]
It is easily seen that \(a_y\) is a \((1,\infty)\)-atom supported in \(P_{y}\). Moreover, 
\[
\begin{aligned}
\sum_{y\in \TT}\mu(y)|f(y)|=\|f\|_{L^1(\mu)}.
\end{aligned}
\]
Therefore \(f\) admits an atomic decomposition and~\eqref{Ieq-H1L1inlem} holds.
\end{proof}

Now, we can turn to the proof of Theorem~\ref{Thm-mainhor}, i.e., giving the $\hardy$-boundedness of $\RR_h$.
\begin{proof}[Proof of Theorem~\ref{Thm-mainhor}]
Let \(f\in \hardy\). By the \(\hardy\to L^1(\mu)\) boundedness of \(\RR_h\) established in Theorem~\ref{Thm-Lpweakttype11q}, we have \(\RR_h f\in L^1(\mu)\). Moreover, by~\eqref{Idn-cancellationq},
\[
\Sigma^*\RR_h f(x)=\frac1q\sum_{y\in s(x)}\RR_h f(y)=0.
\]
Hence, applying Lemma~\ref{Lem-H1L1-cancellation} to \(\RR_h f\), we can deduce that \(\RR_h f\in \hardy\) and combining this with the \(\hardy\to L^1(\mu)\) boundedness of \(\RR_h\) from Theorem~\ref{Thm-Lpweakttype11q}, we obtain
\[
\|\RR_h f\|_{\hardy}\leq \|\RR_h f\|_{L^1(\mu)}\lesssim \|f\|_{\hardy}.
\]
\end{proof}

\section*{Acknowledgements}
The author would like to express her sincere gratitude to Alessio Martini for his helpful advice and insightful discussions during the preparation of this work.

The proof of Theorem \ref{Thm-mainver} was developed with the assistance of ChatGPT-6, and was independently checked by the author, who takes full responsibility for its correctness.

\end{document}